%%% April 2009 %%%%

\input amstex
\input amssym.tex
\documentstyle{amsppt}
\magnification=\magstep1
\hoffset=.25truein
\vsize=8.75truein

\TagsOnLeft

\define\supp{\text{supp}}
\define\Max{\text{Max}}
\define\Min{\text{Min}}
\define\rint{\text{int}}
\define\rloc{\text{loc}}

\define\rVDM{\text{VDM}}
\define\Prob{\text{Prob}}

%\define\rcap{\text{cap}}

\define\cM{{\Cal M}}

\define\IC{{\Bbb C}}
\define\IN{{\Bbb N}}

\define\IR{{\Bbb R}}

\define\prf{\noindent{\bf Proof:}\quad }
%\define\qed{\quad $\square$}

\overfullrule=0pt

\topmatter
\title{Voiculescu's Entropy and Potential Theory}\endtitle
\author{by\\  \\ 
Thomas Bloom*}\endauthor
%\rightheadtext{}
\leftheadtext{}
\thanks{*\ Partially supported by an NSERC of Canada Grant.}
\endthanks

\address{Thomas Bloom, Department of Mathematics, 
University of Toronto, Toronto, Ontario  M5S 3G3 CANADA}
\endaddress
\email{bloom\@math.utoronto.ca}\endemail

%\keywords{}\endkeywords
%\subjclass{}\endsubjclass

%\abstract
%\endabstract
\endtopmatter

\document
\baselineskip=20pt

\head  Introduction\endhead
That the (negative of) the logarithmic energy of
a planar measure can be obtained as a 
limit of volumes
originated with work of D. Voiculescu ([Vo1], [Vo2]).  
His motivation came from
operator theory and free probability theory.
Ben Arous and A. Guionnet [Be-Gu] put that result
in the framework of large deviations. Other results
in that direction are due to Ben Arous and Zeitouni
[Be-Ze] and Hiai and Petz [Hi-Pe]. These authors 
use potential theory and retain the basic form of
Voiculescu's original proof.

Informally, these results express the asymptotic
value (as $d\to \infty$) of the average of a ``weighted"
VanDerMonde determinant of a point
$(\lambda_1, \cdots, \lambda_d)\in E^d$,
as the discrete measures
$\displaystyle \kappa_d (\lambda):={1\over d}\sum^d_{j=1}\delta
(\lambda_j)$ approach a
fixed probability measure $\mu$ with  compact support   $E$ 
in $\IC$. Such weighted VanDerMonde (VDM) determinants
arise, for example, as the joint probability distribution
of the eigenvalues of certain ensembles of random Hermitian matrices
and also in the study of certain determinental point
processes. Specifically, we prove, for measures $\mu$
with support in a rectangle $H$:

\proclaim{Theorem 3.1} 
$$\inf_{G\ni \mu} \lim\limits_{d\to \infty} {1\over d^2}
\log \int_{\tilde G_d(\mu)} | \text{VDM}^w_d (\lambda)^2| d
\tau (\lambda)=\Sigma (\mu) -2 \int Q d \mu$$
where the infimum is over all neighborhoods of $\mu$
in the weak* topology,
$\displaystyle \tilde G_d (\mu):=\{\lambda\in H^d \mid
\kappa_d (\lambda) \in G\},\ Q=-\log w,\ \Sigma (\mu)=
\int\!\!\int\! \log |z-t| d \mu (z) d\mu (t)$,
$\tau$ is a
measure satisfying an appropriate density condition on $H$
(prop 3.1), and $\rVDM^w_d (\lambda)$ is a weighted
VanDerMonde determinant (see (2.7)) with $w$ continuous
and $>0$ on $H$.
\endproclaim

This result is not essentially new however the proof is
new. The lower bound in theorem 3.1 is obtained
by using Markov's polynomial inequality on the
weighted VanDerMondes when the weight is a
real polynomial, the general case being obtained by approximation.

Voiculescu's method (and those of the authors
cited above) uses a ``discretization" argument
on the measure $\mu$ (the method has been used in
other situations ([Ze]-[Ze])). This method relies on
the factorization of the VDM determinent into
linear factors. The method of this paper does not use
such factorization-the interest in doing so, 
being in higher dimensional versions
of these results (The methods of this paper were the
basis for the announcement of some higher dimensional
results [Bl, talk]).

R. Berman ([Be1], [Be2]) has recently proven
large deviation results and a version of the above result
in general higher dimensional
situations. Reduced to the one-dimensional case
of compact subsets of $\IC$, his proof is different
than that of Voiculescu or this paper.
\medskip

\noindent
\centerline{\bf Acknowledgement} 
I would like to thank R. Berman,
L. Bos and N. Levenberg for helpful discussions concerning
this paper.

\head 1.\ Topology on $\cM (E)$ \endhead
Let $E$ be a closed subset of $\IC$ (which we
identify with $\IR^2$). We let $\cM (E)$ denote the set
of positive Borel probability measures on $E$
with the weak* topology.

The weak* topology on $\cM (E)$ is given as follows
(see [E], appendix A8). A neighborhood basis
of any $\mu \in \cM (E)$ is given by sets of the form
$$\{\nu \in \cM (E)\Big| |\int_E  f_i (d\mu - d\nu)| \le \epsilon
\quad {\text{for}}\ \ i=1, \cdots, k\}\tag 1.1$$
where $\epsilon>0$ and $f_1, \cdots, f_k$ are bounded
continuous functions on $E$.

$\cM (E)$ is a complete metrizable space and
for $E$ compact a neighborhood basis of $\mu \in \cM (E)$
is given by sets of the form
$$G(\mu, k, \epsilon):=\{\nu \in \cM (E)\Big|
| \int_E x^{n_1} y^{n_2} (d\mu - d\nu)| <\epsilon\}
\quad {\text{for}}\ \ k, n_1, n_2 \in \IN,\ 
n_1 + n_2 \le k\quad {\text{and}}\ \ \epsilon >0.\tag 1.2$$

That is $G(\mu, k, \epsilon)$ consists of all probability
measures on $E$ whose (real) moments, up to order $k$,
are within $\epsilon$ of the corresponding moment for $\mu$.

It is clear that for $k_1 \ge k$ and $\epsilon_1 \le \epsilon$ that
$$G(\mu, k_1, \epsilon_1) \subset G(\mu, k, \epsilon).\tag 1.3$$
Now for $\lambda=(\lambda_1, \cdots, \lambda_d)\in \IC^d$, we let
$$\kappa_d (\lambda):={1\over d} \sum^d_{j=1} \delta (\lambda_j)
\tag 1.4$$
where $\delta$ is the Dirac $\delta$-measure at the
indicated point.

We let
$$\tilde G_d(\mu, k, \epsilon):=\{\lambda\in E^d\Big|
\kappa_d (\lambda) \in G(\mu, k, \epsilon)\}.\tag 1.5$$

It follows from (1.3) that
$$\tilde G_d(\mu, k_1, \epsilon_1)\subset \tilde G_d
(\mu, k, \epsilon)\quad {\text{for}}\ \ k_1\ge k \quad 
{\text{and}}\ \ \epsilon_1\le \epsilon.\tag 1.6$$
For $\lambda \in \IC^d$ we let
$$\Delta_d (\lambda)=\{\lambda' \in \IC^d \Big| |\lambda'_j -
\lambda_j| \le e^{-\sqrt d}\quad {\text{for}}\ \ j=i, \cdots, d\}.
\tag 1.7$$

Proposition 1.1 and 1.2 follow immediately 
from the definition of the weak* topology in $\cM (E)$
(for $E$ compact).

\proclaim{Proposition 1.1}
Let $f$ be continuous on $E$ and $\mu \in \cM (E)$.
Given $\epsilon_1 > 0$ there exist $k, \epsilon$ such that
$$\Big| \int_E f(d\mu - \kappa_d (\lambda))\Big| \le \epsilon_1
\quad {\text{for}}\ \ \lambda\in \tilde G_d(\mu, k, \epsilon).$$
\endproclaim

\proclaim{Proposition 1.2}
Let $\nu \in G(\mu, k, \epsilon)$. Then there
exists $k_1, \epsilon_1$, such that
$G(\nu, k_1, \epsilon_1) \subset G(\mu, k, \epsilon)$.
\endproclaim

\proclaim{Proposition 1.3}
Let $\lambda \in \tilde G_d (\mu, k, \epsilon)$. Then
$\Delta_d (\lambda)\in \tilde G_d (\mu, k, 2 \epsilon)$
for all $d$ sufficiently large.
\endproclaim

\prf The proof follows from the fact that
monomials  satisfy a Lipshitz condition on $E$.

\head 2.\ Markov's Polynomial Inequality  \endhead
The classical Markov polynomial inequality for
real polynomials on an interval $I\subset \IR$ is an 
estimate for the derivative of the polynomial
in terms of its degree and sup norm on $I$.
Specifically ([Be-Er], theorem 5.1.8)
$$| p'(x) | \le A k^2 \| p\|_I\quad {\text{for}}\ \ 
x\in I\tag 2.1$$
where $k=\deg (p)$ and $A$ is a constant $>0$. For
$I=[-1, 1]$ on may take $A=1$.

Numerous extensions of (2.1) to multivariable
settings have been established (see e.g. [Ba], [Pl]).

We will however use a version of (2.1)
for rectangles $H\subset \IR^2$ which is an immediate
consequence of (2.1). (We will always assume that
rectangles have sides parallel to the axes). Let
$p(x, y)$ be a polynomial of degree $\le k$ in each varaible, then
$$| {\text{grad}}\  (p)(x)| \le Ak^2 \| p\|_H\tag 2.2$$
where $A>0$ is a constant.

Integrating (2.2) over the straight line
joining $z_1$ to $z_2$ in $H$ we have
$$| p(z_1)-p(z_2)| \le Ak^2 \|p\|_H |z_1-z_2|.\tag 2.3$$

We will now use (2.3) to show in
quantitative terms that the value of polynomials
at points near a point where it assumes its
maximum is close to the maximum value.

Let $\{\Lambda_d\}_{d=1,2\cdots}$ be a sequence of
polynomials on $(\IR^2)^d$, non negative on $H^d$,
such that for some constants $c_1>0,\ \gamma_1>0$ each
polynomial $\Lambda_d$ is of degree $\le c_1 d^{\gamma_1}$ in 
each of its $2d$ real variables. Let
$z^M:=(z^M_1, \cdots, z^M_d)$ be a point in
$H^d\subset \IC^d \simeq \IR^{2d}$ where $\Lambda_d$
assumes its maximum
i.e. $\Lambda_d (z^M)=\| \Lambda_d\|_{H^d}$.

\proclaim{Theorem 2.1}
For $z\in \Delta_d (z^M)\cap H^d$. Then
$$\Lambda_d (z) \ge \Lambda_d (z^M) \psi (d)$$
where $\psi (d)=1-c d^\gamma e^{-\sqrt d}$ for some
constants $c, \gamma >0$ (independent of $d$).
\endproclaim

\prf We write $\Lambda_d (z^M) -\Lambda_d (z)$ in the form
$$\Lambda_d(z^M)-\Lambda_d(z)=\sum^d_{j=1} \Lambda_d 
(z_1, \cdots, z_{j-1}, z^M_j, \cdots, z^M_d)-\Lambda_d
(z_1, \cdots, z_j, z^M_{j+1}, \cdots, z^M_d).\tag 2.4$$
But for $z_1, \cdots, z_{j-1}, z^M_{j+1}, \cdots, z^M_d$ fixed,
$t\to \Lambda_d (z_1, \cdots, z_{j-1}, t, z^M_{j+1}, \cdots,
z^M_d)$ is a polynomial
in $\IR^2$ of $\deg\le c_1 d^{\gamma_1}$ in each real variable.
Applying (2.3) and the fact that $z\in \Delta_d (z^M)$ to each
term on the right side of (2.4) we have an estimate of the form
$$\Lambda_d (z^M) -\Lambda_d (z)\le dA(c_1 d^{\gamma_1})^2 \Lambda_d
(z^M) e^{-\sqrt d}.\tag 2.5$$
The result follows. 
\qed

We will apply this result to sequences of polynomials 
constructed as follows: Let
$$\rVDM_d (\lambda)=\rVDM_d (\lambda_1, \cdots, \lambda_d)=
\prod_{1\le i<j\le d} (\lambda_i-\lambda_j)\tag 2.6$$
be the VanDerMonde determinant of the points
$(\lambda_1, \cdots, \lambda_d)$ in $\IC$.
Also, for $w$ a function on $\IC$ we let
$$\rVDM^w_d (\lambda):=\rVDM_d (\lambda)
{\displaystyle{d\atop\prod}\atop{\scriptstyle i=1}} 
w (\lambda_i)^d.\tag 2.7$$
Thus if $w$ is a real polynomial the sequence of 
polynomials $\Lambda_d:=|\rVDM^w_d (\lambda)|^2$ for
$d=1, 2, \cdots$ satisfies the hypothesis of theorem 2.1.
In this situation, points $z^M \in H^d$ at which 
$| \rVDM^w_d (\lambda)|^2$ assumes its maximum are known as
a $w$-Fekete set.

\head 3.\ Energy as a limit of volumes  \endhead
Let $E$ be a compact subset of $\IC$ and 
$w$ an admissible weight function on $E$ (i.e. $w$
is uppersemicontinuous, $w\ge 0,\ w>0$ on a non-polar 
subset of $E$.  In particular, $E$ is non-polar.

The weighted equilibrium measure (see [Sa-To], theorem I 1.3),
denoted $\mu_{eq} (E, w)$ is the unique probability measure
which minumizes the functional $I_w (\nu)$ over all
$\nu \in \cM (E)$
where
$$\eqalign{I_w(\nu):& =\int\!\!\int \log
\Big({1\over |z-t|w(z)w (t)}\Big) d\nu (z) d\nu (t)\cr
   & =-\int\!\!\int \log |z-t| d\nu (z) d\nu(t)+2\int Q(z)d\nu(z)\cr}
\tag 3.1$$
where
$$ Q(z):=-\log w (z). \tag 3.2$$

$I_w (\nu)$ is termed the weighted energy of the measure 
$\nu$.  We also use the notation
$$\Sigma (\nu):=\int\!\!\int \log |z-t| d\nu (z) d\nu (t). \tag 3.3$$

$\Sigma (\nu)$ is termed the free entropy of $\nu$ (it may assume the
value $-\infty$). We let
$$\delta^w_d:=\Max_{\lambda\in E^d} | \rVDM^w_d
(\lambda)|^{2\over d(d-1)}. \tag 3.4$$
Then (see [Sa-To], chapter III, theorem 1.1)
$$\delta^w:=\lim\limits_{d\to \infty} \delta^w_d \tag 3.5$$
exists and
$$\log \delta^w=-I_w (\mu_{eq} (E, w))=\Sigma (\mu_{eq} (E, w)) -
2\int Q(z)d\mu_{eq} (E, w).\tag 3.6$$

Now, let $\tau$  be a positive Borel measure on $E$.

We say that the triple $(E, w, \tau)$
satisfies the weighted Bernstein-Markov (B-M)
inequality if, for all $\epsilon >0$, these exists a constant
$c> 0$ such that, for all analytic polynomials
$p$ of degree $\le k$ we have
$$\| w^k p\|_E \le c(1+\epsilon)^k \| w^k p\|_{L^2 (\tau)}.\tag 3.7$$

We set
$$Z_d:=\int_{E^d} | \rVDM^w_d (\lambda)|^2 d\tau(\lambda).\tag 3.8$$
where $d\tau (\lambda)=d\tau (\lambda_1) \cdots d\tau(\lambda_d)$
is the product
measure on $E^d$. Then if $(E, w, \tau)$ satisfies the
weighted B-M inequality ([Bl-Le2]).
$$\lim\limits_{d\to \infty} Z^{d^{-2}}_d =\delta^w.\tag 3.9$$

We will need measures on $E$ which satisfy
the weighted B-M ineqality for all continuous
admissible weights.

To this end we consider measures $\tau$ 
which satisfy the following condition (satisfied by any measure
that is a positive continuous function times Lebesgue measure):

There is a constant $T>0$ such that
$$\tau (D(z_0, r)) \ge r^T\quad {\text{for\ all}}\ \
z_0 \in E\quad {\text{and}}\ \ r\le r_0. \tag 3.10$$

Here $D(z_0, r)$ denotes the disc center $z_0$ radius $r$ and
$r_0 >0$.

\proclaim{Proposition 3.1}
Let $H$ be a rectangle in $\IC$ and let
$\tau$ satisfy (3.10). Then for all continuous functions
$w>0$ on $H,\ (H, w, \tau)$ satisfies the weighted B-M
inequality.
\endproclaim

\prf First we can consider $H\subset \IC \simeq \IR^2$ as
a subset of $\IC^2$.
Then, using [Bl-Le1], theorem 2.2 and [B1], theorem 3.2.
$(H, w, \tau)$ satisfies the weighted B-M inequality as a
subset of $\IC^2$ (the definition of which is an obvious adaptation of (3.7)
to the several variable case-see [Bl]). But every analytic
polynomial $p(z)$ on $\IC$ is the restriction to
$\IC\simeq \IR^2 \subset \IC^2$ of the analytic
polynomial $p(z_1 + i z_2)$. Hence the result \qed

Let $H$ be a rectangle in $\IC$, $\mu \in \cM (H)$, $\tau$ satisfy
(3.10), and let $\phi>0$ be a continuous function on $H$.
Let $S=-\log \phi$. We will consider integrals
of the form
$$J^\phi_d (\mu, k, \epsilon):= \int_{\tilde G_d (\mu, k, \epsilon)}
| \rVDM^\phi_d (\lambda)|^2 d\tau (\lambda). \tag 3.11$$

The integral in (3.11) is of the same form as
that in (3.8) used to define $Z_d$ however here
we only integrate over a subset of $H^d$. Theorem
3.1 below establishes asymptotic properties of
such integrals. The leading term depends only on $\mu$ on $S$
(and as mentioned in te introduction, the result is not
essentially new
but goes back to results of Voiculescu ([Vo1], [Vo2]).

\proclaim{Theorem 3.1}
$$\inf_{k, \epsilon}\Big\{\lim\limits_{d\to\infty} {1\over d^2}
\log J^\phi_d (\mu, k, \epsilon)\Big\}=\Sigma (\mu)-
2\int S d\mu$$
\endproclaim

\prf To prove this result we will show
\smallskip
\item{(a)\ }
$\displaystyle \inf_{k, \epsilon}
\Big\{{\varlimsup_{d\to\infty}}{1\over d^2} \log
J^\phi_d (\mu, k, \epsilon)\Big\} \le \Sigma
(\mu)-2\int Sd\mu$\qquad  and
\medskip

\item{(b)\ }
$\displaystyle \inf_{k, \epsilon}
\Big\{{\varliminf_{d\to\infty}}{1\over d^2} \log
J^\phi_d (\mu, k, \epsilon)\Big\} \ge \Sigma (\mu)-2\int Sd\mu.$
\medskip

To prove the upper bound (a) we will
follow ([Be1], proposition 3.4). The proof does not use (3.10).
Let $w$ be continuous $>0$ on $H$.
Then
$${\displaystyle{d\atop\prod}\atop{\scriptstyle i=1}} 
w(\lambda_i)^{2d} |\rVDM^\phi_d (\lambda)|^2
= |\rVDM^w_d (\lambda)|^2 
{\displaystyle{d\atop\prod}\atop{\scriptstyle i=1}} 
\phi (\lambda_i)^{2d}.
\tag 3.12$$
Hence,
$$| \rVDM^\phi_d (\lambda)|^2 \le (\delta^w_d)^{d(d-1)}
\exp \big(2d^2\int_H (Q-S) \kappa_d (\lambda)\big).\tag 3.13$$

Let $\lambda^d \in \overline{{\tilde G_d} (\mu, k, \epsilon)}$ be
a point at which the 
maximum of $| \rVDM^\phi_d (\lambda)|$ over
$\overline{\tilde G_d (\mu, k, \epsilon)}$ is attained.
(3.13) implies that
$$J^\phi_d (\mu, k,\epsilon) \tau (H)^d \le(\delta^w_d)^{d(d-1)}
\exp (2d \int_H (Q-S) \kappa_d (\lambda^d)).\tag 3.14$$

For any sequence of $d$'s we may pass to a
subsequence and assume that the sequence of
measures $\kappa_d (\lambda^d)$ converges to a measure
$\sigma \in \overline{G(\mu, k, \epsilon)}$.

We deduce that
$${\varlimsup_{d\to \infty}} {1\over d^2} \log
J^\phi_d (\mu, k, \epsilon)\le \log \delta^w +2\int_H
(Q-S) d\sigma.\tag 3.15$$

Taking the $\inf$ over $k, \epsilon$, the $\sigma$'s converge
to $\mu$ so
$$\inf_{k, \epsilon} {\varlimsup} {1\over d^2}\log J^\phi_d
(\mu, k, \epsilon)\le \log \delta^w +2 \int_H (Q-S) d\mu.$$

Now take a sequence of continuous weights $w$ such that
$\mu_{eq} (H, w)$ converges to $\mu$ in $\cM (H)$ and
$\Sigma (\mu_{eq} (H, w))$ converges to $\Sigma (\mu)$
(see proof of (b) (iii)).

Then using (3.6) we obtain (a).

For the lower bound (b) we proceed as follows.

We prove (b) when

\item{(i)\ }
$\mu=\mu_{eq} (H, w)$, $w$ is a polynomial $>0$ and $\phi=w$.

\item{(ii)\ }
$\mu$ as in (i) but the restriction on $\phi$ is dropped.

\item{(iii)\ }
general $\mu$.

\item{(i)\ }
We consider points $z^M \in H^d$ at which
$| \rVDM^w_d (\lambda)|^2$ assumes its maximum 
(i.e. $w$-Fekete points).  It is known
that $\kappa_d (z^M)$ converges to $\mu$ in $\cM(H)$ so, for $d$
large, $\kappa_d (z^M)\in \tilde G_d (\mu, k, \epsilon)$. Then
for $d$ sufficiently large, using proposition 1.3
$$J^w_d (\mu, k, 2\epsilon) \ge \tau (\Delta_d(z^M)) 
\Min_{\lambda\in \Delta_d(z^M)\cap H} | \rVDM^w_d (\lambda) |^2.
\tag 3.16$$

By (3.10) $\tau(\Delta_d (z^M))\ge e^{-Td\sqrt d}$ and using
theorem 2.1 on the sequence of polynomials $| \rVDM^w_d (\lambda)|^2$
we have
$$\varliminf_{d\to\infty} {1\over d^2} \log J^w_d
(\mu, k, 2\epsilon) \ge \log \delta^w =\Sigma (\mu)-2 
\int Q d\mu.\tag 3.17$$

\item{(ii)\ }
Given $\epsilon_1 >0$, by proposition 1.1, choose $k, \epsilon$
so that 
$$\int (Q-S) (d\mu - \kappa_d (\lambda)) \le\epsilon,\quad
{\text{for\ all}} \ \ \lambda\in \tilde G_d (\mu, k, \epsilon).
\tag 3.18$$

This yields
$${\displaystyle{d\atop\prod}\atop{\scriptstyle i=1}} 
w (\lambda_i)^{2d} \le
{\displaystyle{d\atop\prod}\atop{\scriptstyle i=1}} 
\phi (\lambda_i)^{2d} \exp (2d^2
[\epsilon_1 -\int (Q-S) d\mu]).\tag 3.19$$

Multiplying by $| \rVDM_d (\lambda)|^2$ and integrating over
$\tilde G_d (\mu, k, \epsilon)$ 
gives
$$\exp (2d^2 [-\epsilon_1 +\int (Q-S) d\mu]) J^w_d (\mu, k, \epsilon)
\le J^\phi_d (\mu, k, \epsilon).\tag 3.20$$

Then using (i) and the fact that $\epsilon_1 >0$ is arbitrary gives
$$\inf_{k, \epsilon} \varliminf_{d\to \infty} 
\Big\{{1\over d^2} \log J^\phi_d (\mu, k, \epsilon) \Big\}
\ge \log \delta^w +2\int (Q-S) d\mu\tag 3.21$$
and using (3.6) completes (ii).

For (iii) we will use an approximation argument.
First we note that it is an immediate consequence of
proposition 1.2 that
$$\mu\to \inf_{k, \epsilon} \varliminf_{d\to \infty} {1\over d^2}
\log J^\phi_d (\mu, k, \epsilon)\quad {\text{is}}$$
uppersemicontinuous on $\cM (H)$. So, it suffices to show
that any $\mu\in \cM (H)$ may be approximated by
measures $\{\mu_s\}\in \cM (H)$ where each $\mu_s$ satisfies 
(i) above and $\Sigma(\mu_s)$ converges to $\Sigma (\mu)$. 
$(\mu\to \Sigma (\mu)$ is uppersemicontinuous
on $\cM (H)$ but not, in general, continuous). First, we may
assume $\supp (\mu) \subset \rint (H)$ since, taking $H$
centered at $0$ the measures $\mu_s =\pi^*_s (\mu)$, the
push forward of $\mu$ under the scaling $z\to sz (s<1)$
satisfy $\Sigma(\mu_s)$
converges to $\Sigma (\mu)$. Next for $\mu$ with compact support in
$\rint (H)$ we approximate $\mu$ by $\mu_s=\mu*\rho_s$
where $\rho=\rho(|z|)$ is
a standard smoothing kernel for subhamonic functions
on $\IC$ and $\rho_s=s^{-2} \rho\Big({|z| \over s}\Big)$.
Then $\rho(|z|)\rho(|t|)$ is
a standard smoothing kernel for plurisubhamonic functions
on $\IC^2$ so $\log |z-t| * \Big(\rho_z (|z|) \rho_s(|t|)\Big)$
decreases pointwise to $\log|z-t|$. 

Now, for $\nu$ a positive measure with compact 
support in $\IR^n$, $\psi$ a smooth function with
compact support such that $\psi(x)=\psi(-x)$ and 
$h\in L^1_{\rloc} (\IR^n)$ then
$$\int_{\IR^n} \psi (\nu * h) dm \ =\int_{\IR^n} (\psi * h) d\nu$$
where $dm$ denotes Lebesque measure and $*$ convolution.

Applying this formula to $\IR^4\simeq \IC^2$ with $(z, t)$ as
coordinates, $\nu=\mu \otimes \nu$, $\psi=\rho_s (|z|) \rho_z (|t|)$ 
and $h=\log |z-t|$, then
using the Lebesgue monotone convergence theorem yields
$\Sigma (\mu_s)\to \Sigma(\mu)$ (as $s\to 0)$.
Finally, for $\mu$ a smooth function with compact support
times Lebesgue measure let $Q$ be a smooth potential for $\mu$
which, adding a constant, we may assume is $<0$ on $H$.
Then $\mu=\mu_{eq} (H, w)$ where $w=e^{-Q}$ and one may
approximate $\mu$ by $\mu_s=\mu_{eq} (H, w_s)$ where $w_s$ are real
polynomial weights converging uniformly to $w$ on $H$.  To see that
$\Sigma(\mu_s)$ converges to $\Sigma(\mu)$ we may use
([Sa-To], theorem 6.2 (c), chapter I) - which is
stated for monotonically decreasing 
sequences of weights but the conclusion also
holds for uniformly convergent sequences of weights.

\head 4.\ Entropy  \endhead
Let $\mu \in \cM (H)$. The free entropy of $\mu$ (see (3.3))
defined as an integral may be obtained via discrete
measures as follows:

Let $W(\mu)$ be defined via
$$W_d (\mu, k,\epsilon):=\sup 
\{| \rVDM_d(\lambda)|{^{2\over d(d-1)}}\Big|
\kappa_d (\lambda)\in \tilde G_d (\mu, k, \epsilon)\}\tag 4.1$$
and let
$$W (\mu, k, \epsilon)={\varlimsup_{d\to \infty}} W_d
(\mu, k,\epsilon)\tag 4.2$$
and
$$ W(\mu)=\inf_{k,\epsilon} W(\mu, k, \epsilon).\tag 4.3$$
Then

\proclaim{Theorem 4.1}\ $\log W(\mu)=\Sigma (\mu)$.
\endproclaim

\prf The proof consists of establishing the 
two inequalities
\item{(a)} $\log W(\nu) \le \Sigma(\mu)$ and
\item{(b)} $\Sigma(\mu)\le \log W(\mu).$

For (a) let $\kappa_d (\lambda^d) ={1\over d} \sum^d_{j=1}
\delta (\lambda^d_j)$ be, for $d=1,2, \cdots ,$ a sequence of
discrete measures converging to $\mu$
 weak* such that
$$\log W(\mu)=\lim\limits_{d\to \infty} {1\over d^2}\sum_{j\not= k}
\log | \lambda^d_j - \lambda^d_k|.\tag 4.4$$
Now,
$$\lim\limits_{d\to \infty} {1\over d^2} \sum_{j\not= k} \delta
(\lambda^d_j, \lambda^d_k)=\mu\otimes\mu \ {\text{weak}}^*\tag4.5$$
and so, since $\log|z-t|$ is u.s.c.
\item{(a)} follows from ([Sa-To], theorem 1.4, chapter O).

For (b) we note that by definition of the quatities involved
$$J_d (\mu, k, \epsilon)\leq W_d (\mu, k, \epsilon)^{d(d-1)} 
{\tau}(H)^d.$$
so that
$$\varlimsup_{d\to\infty} {1\over d^2}\log J_d
(\mu, k, \epsilon)\le \log W(\mu, k, \epsilon)$$
Taking the $\inf$ over $k, \epsilon$ and using
theorem 3.1, (b) follows.\qed

\proclaim{Corollary 4.1}
Define $W_d^{\phi} (\mu, k, \epsilon)$ analogously to 
the definition of $W_d (\mu, k, \epsilon)$ in (4.1). 
That is
$W_d^{\phi} (\mu, k, \epsilon) =\sup \{| \rVDM_d^{\phi}
(\lambda)|^{2\over d(d-1)}| \kappa_d (\lambda) \in \tilde G_d
(\mu, k, \epsilon)\}$ and define $W^\phi (\mu)$ analogously to the
definition of $W(\mu)$ (see (4.3)). Then
$$\Sigma (\mu)-2\int Sd\mu=\log W^\phi (\mu).$$
\endproclaim

\head 5.\ Large Deviation  \endhead
Consider the sequence of probability measures
on $H^d$ (for $d=1, 2, \cdots)$ given by
$${| \rVDM^\phi_d (\lambda)|^2 d\tau (\lambda)\over
       Z^\phi_d} :=\Prob_d.\tag 5.1$$
       Then
$${1\over d^2} \log \Prob_d (\tilde G_d (\mu, k, \epsilon))=
\Prob_d \{\lambda | \kappa_d (\lambda)\in G(\mu, k, \epsilon)\}.
\tag 5.2$$

Using theorem 3.1 and (3.9) gives
$$\inf_{k,\epsilon}\lim\limits_{d\to\infty} 
{1\over d^2} \log \Prob_d
(\tilde G_d (\mu, k, \epsilon)) =I_\phi (\mu_{eq} (H, \phi))
- I_\phi (\mu).$$
The functional $\mu\to I_\phi (\mu)-I_\phi (\mu_{eq} (H,\phi))=:
I(\mu)$
attains its minimum value of zero at the unique
measure $\mu=\mu_{eq} (H, \phi)$.

Then $I(\mu)$ is a good rate functional and the sequence of
descrete random measures $\kappa_d(\lambda)$ satisfy a large
devation principle in the scale $d^{-2}$ (see discussion 
[Hi-Pe], page 211).
\vfill\eject

\vglue .1truein

Thomas Bloom

Department of Mathematics

University of Toronto

Toronto, ON

CANADA\ \ M5S 2E4

email: bloom\@math.toronto.edu

\end
\Refs
\widestnumber\key{ABCDE}
\smallskip

\ref\key Ba 
\by M. Baran 
\paper Complex equilibrium measure and Bernstein type theorems for
compact sets in $\IR^n$
%\inbook Recent Developments in Several Complex Variables
\jour Proc Am. Math. Soc.   
\vol 123. no. 2 
\yr 1995 
\pages 485-494 
\endref

\ref\key Be1 
\by R. Berman 
\paper Large deviations and entropy for determinental point processes
on complex manifolds. 
\vol arxiv:0812.4224
\endref

\ref\key Be2 
\by R. Berman 
\paper Determinental point processes and fermions on complex
manifolds: bulk university
\vol arxiv:0811.3341
\endref


\ref\key B1-Le1 
\by T. Bloom and N. Levenberg 
\paper Capacity convergence results and applications to
a Bernstein-Markov inequality
%\inbook Recent Developments in Several Complex Variables
\jour Tr. Am. Math Soc.   
\vol 351. no. 12
\yr 1999 
\pages 4753-4767
\endref

\ref\key B1-Le2 
\by T. Bloom and N. Levenberg 
\paper Asymptotics for Christoffel functions of Planar Measures 
\jour J. D'Anal Math.
\vol 106
\yr 2008
\pages 353--371
\endref

\ref\key B1-Le3 
\by T. Bloom and N. Levenberg 
\paper Transfinite diameter notions in $\IC^n$ and integrals
of VanDerMonde determinants 
\vol arxiv:0712.2844
\endref

\ref\key B1 
\by T. Bloom  
\paper Weighted polynomials and weighted pluripotential theory
%\inbook Recent Developments in Several Complex Variables
\jour Tr. Am. Math Soc.   
\vol 361 no. 4 
\yr 2009 
\pages 2163-2179 
\endref

\ref\key B1, talk 
\by T. Bloom 
%\paper
\inbook "Large Deviations for VanDerMonde determinants" talk given at the
Workshop on Complex Hyperbolic Geometry and Related Topics (November 17-21,
2008) at the Fields Institute. 
http://www.fields.utoronto.ca/audio/08-09/hyperbolic/bloom/
\endref

\ref\key Be-Gu 
\by G. Ben Arous and A. Guionnet
\paper Large deviation for Wigner's law and Voiculescu's
non-commutative entopy
\jour  Prob. Th. Related Fields
\vol 108 
\yr 1997
\pages 517-542 
\endref

\ref\key Be-Ze 
\by G. Ben Arous and O. Zeitouni 
\paper Large deviations from the circular law
\jour ESAIM: Probability and Statistics 
\vol 2 
\yr 1998 
\pages 123-134 
\endref


\ref\key Bo-Er 
\by P. Borwein and T. Erdelyi 
\paper Polynomials and Polynomial Inequalities 
\jour Springer 
Graduate Texts in Mathematics 
\vol 161 
\yr 1995 
\publaddr New York
\endref

\ref\key De-Ze
\by A. Dembo and O. Zeitouni 
\paper Large Deviation Techniques and Applications, 2nd edition
\jour  Springer 
\yr 1998
\publaddr New York
\endref

\ref\key Dei
\by P. Deift
\paper Othogonal Polynomials and Random Matrices: A Riemann-Hilbert
approach
\jour AMS Providence RI
\yr 1999
\endref

\ref\key El
\by R.S. Ellis
\paper Entropy, Large Deviations and Statistical Mechanics
\jour Springer 
\publaddr New York/Berlin (1985)
\endref

\ref\key Hi-Pe
\by F. Hiai and D. Petz
\paper The Semicircle Law, Free Random Variables and Entropy
\jour AMS Providence RI (2000)
\endref

\ref\key Kl
\by M. Klimek
\paper Pluripotential Theory
\jour Oxford University Press
\yr 1991
\publaddr Oxford
\endref

\ref\key Pl
\by W. Plesniak
\paper Inegalite de Markov en plusieurs variables
\jour Int. J of Math and Math. Science, Art ID 24549 
\yr 2006
\pages 1--12
\endref

\ref\key Sa-To
\by E. Saff and V. Totik
\paper Logarithmic Potential with External Fields
\jour Springer-Verlag
\publaddr Berlin (1997)
\endref

\ref\key St-To
\by H. Stahl and V. Totik
\paper General Orthogonal Polynomials
\jour Cambridge Unviersity Press  
\publaddr Cambridge (1992)
\endref

\ref\key Vo1
\by D. Voiculescu
\paper The analogues of entropy and of Fisher's information measure
in free probability theory I
\jour Comm. Math. Phy. 
\vol 155
\yr 1993
\pages 71-92
\endref

\ref\key Vo2
\by D. Voiculescu
\paper The analogues of entropy and of Fisher's information
measure in free probability II
\jour Inv. Math.
\vol 118
\yr 1994
\pages 411--440
\endref

\ref\key Ze-Ze
\by O. Zeitouni and S. Zelditch
\paper Large Deviations of empirical zero
point measures on Riemann surfaces, $I:g=0$ 
\jour arxiv:0904, 4271 
\endref
